\documentclass[12pt,abstract]{scrartcl}
\usepackage[dvipdfmx]{graphicx}
\usepackage{amsmath,amsthm}
\usepackage[round]{natbib}
\newtheorem{proposition}{Proposition}
\usepackage{newtxtext}
\usepackage{newtxmath}
\begin{document}
\title{Dynamic monopolistic competition with sluggish adjustment of entry and exit}
\author{Yasuhito Tanaka\\[.2cm]
Faculty of Economics, Doshisha University,\\
 Kamigyo-ku, Kyoto, 602-8580, Japan.\\[.1cm]
E-mail:yatanaka@mail.doshisha.ac.jp}

\date{}
\maketitle

\begin{abstract}
We study a steady state of a free entry oligopoly with differentiated goods, that is, a monopolistic competition, with sluggish adjustment of entry and exit of firms under general demand and cost functions by a differential game approach. Mainly we show that the number of firms at the steady state in the open-loop solution of monopolistic competition is smaller than that at the static equilibrium of monopolistic competition, and that the number of firms at the steady state of the memoryless closed-loop monopolistic competition is larger than that at the steady state of the open-loop monopolistic competition, and may be larger than the number of firms at the static equilibrium.
\end{abstract}
\begin{description}
\item[Keywords:] monopolistic competition; differential game; general demand function; general cost function; open-loop; closed-loop
\end{description}

\section{Introduction}

There are many studies of an oligopoly by differential game theory, for example, \cite{cl0}, \cite{cl6}, \cite{cl1}, \cite{cl3}, \cite{cl2}, \cite{cl4}, \cite{fu}, \cite{fu1} and \cite{lam18}. Most of these studies used a model of specific (linear or exponential) demand functions and specific (quadratic or linear) cost functions. We study a steady state of a dynamic free entry oligopoly with differentiated goods, that is, a monopolistic competition with sluggish adjustment of entry and exit of firms under general demand and cost functions by a differential game approach. In the next section we present a model and assumptions. We consider a dynamics of the number of firms which enter into the industry according to the rule that the number of firms increases or decreases proportionally to the total profits of the firms\footnote{Alternatively, we can assume that the number of firms increases or decreases proportionally to the average profit of the firms. Essentially the same result is obtained in both cases.}. In Section 3 we consider an open-loop solution of a differential game analysis of monopolistic competition. We present both a general analysis and a linear example. In Section 4 we examine a general model of a memoryless closed-loop solution. In Section 5 we consider an example with linear demand and cost functions of the memoryless closed-loop solution. We compare open-loop and memoryless closed-loop solutions, and mainly show the following results.  
\begin{enumerate}
	\item The number of firms at the steady state in the open-loop solution of monopolistic competition is smaller than that at the static equilibrium of monopolistic competition.
	\item The number of firms at the steady state in the memoryless closed-loop solution of monopolistic competition is larger than that at the steady state of the open-loop solution of monopolistic competition, and may be larger than the number of firms at the static equilibrium.
\end{enumerate}

We also show that when the discount rate (denoted by $\rho$) approaches to positive infinity, or the speed of adjustment of the number of firms approaches to zero, the steady states of the open-loop and the closed-loop solutions approach to the static equilibrium of monopolistic competition.

\section{The model and free entry condition}

There is a symmetric oligopoly where, at any $t\in [0, \infty)$, $n$ firms, Firms 1, 2, $\dots$, $n$ produce differentiated goods. The firms maximize their discounted profits. Let $x_i(t),\ i\in \{1, 2, \dots, n\}$, be the outputs of the firms, $p_i(t)$ be the price of the good of Firm $i$ at $t$. 

The inverse demand function for Firm $i,\ i\in \{1, 2, \dots, n\}$, is
\[p_i(t)=p_i(x_1(t), x_2(t), \dots, x_n(t)),\ i\in \{1, 2, \dots, n\}.\]
For simplicity we denote $p_i(x_1(t), x_2(t), \dots, x_n(t))$, $\frac{\partial p_i(x_1(t), x_2(t), \dots, x_n(t))}{\partial x_i(t)}$, $\frac{\partial p_i(x_1(t), x_2(t), \dots, x_n(t))}{\partial x_j(t)}$,\\
$\frac{\partial^2 p_i(x_1(t), x_2(t), \dots, x_n(t))}{\partial x_i(t)^2}$, $\frac{\partial^2 p_i(x_1(t), x_2(t), \dots, x_n(t))}{\partial x_i(t)\partial x_j(t)},\ j\neq i$, by $p_i$, $\frac{\partial p_i}{\partial x_i(t)}$, $\frac{\partial p_i}{\partial x_j(t)}$, $\frac{\partial^2 p_i}{\partial x_i(t)^2}$, $\frac{\partial^2 p_i}{\partial x_i(t)\partial x_j(t)}$, and so on.
We assume
\[\frac{\partial p_i}{\partial x_i(t)}<0,\ i\in \{1, 2, \dots, n\},\]
\[\frac{\partial p_i}{\partial x_j(t)}<0,\ j\neq i,\]
\[\left|\frac{\partial p_i}{\partial x_j(t)}\right|<\left|\frac{\partial p_i}{\partial x_i(t)}\right|,\]
and
\begin{equation}
\frac{\partial p_i}{\partial x_j(t)}+\frac{\partial^2 p_i}{\partial x_i(t)\partial x_j(t)}x_i(t)<0,\ i\in \{1, 2, \dots, n\},\ j\neq i.\label{sub1}
\end{equation}
The last condition means that the outputs of the firms are strategic substitutes. Note that
\[\frac{\partial p_i}{\partial x_j(t)}+\frac{\partial^2 p_i}{\partial x_i(t)\partial x_j(t)}x_i(t)=\frac{\partial^2 p_ix_i(t)}{\partial x_i(t)\partial x_j(t)}.\]
Similarly,
\[\frac{\partial p_i}{\partial x_j(t)}+\frac{\partial^2 p_i}{\partial x_j(t)\partial x_k(t)}x_k(t)=\frac{\partial^2 p_ix_k(t)}{\partial x_j(t)\partial x_k(t)},\ j\neq i,\ k\neq i, j.\]
We assume
\[\left|\frac{\partial^2 p_ix_i(t)}{\partial x_i(t)\partial x_j(t)}\right|\geq \left|\frac{\partial^2 p_ix_k(t)}{\partial x_j(t)\partial x_k(t)}\right|.\]
Then, we obtain
\begin{equation}
\frac{\partial p_i}{\partial x_j(t)}+\frac{\partial^2 p_i}{\partial x_j(t)\partial x_k(t)}x_k(t)<0,\ j\neq i,\ k\neq i, j.\label{sub2}
\end{equation}
and
\begin{equation}
\left|\frac{\partial p_i}{\partial x_j(t)}+\frac{\partial^2 p_i}{\partial x_i(t)\partial x_j(t)}x_i(t)\right|\geq \left|\frac{\partial p_i}{\partial x_j(t)}+\frac{\partial^2 p_i}{\partial x_j(t)\partial x_k(t)}x_k(t)\right|.\label{sub3}
\end{equation}
By symmetry of the model at the steady states of open-loop and closed-loop solutions $x_i(t)=x_j(t)=x_k(t)$.

About the derivative of $p_i$ with respect to $n$ we have
\[\frac{\partial p_i}{\partial n(t)}=\frac{\partial p_i}{\partial x_j(t)}x_j(t).\]
The cost function of Firm $i,\ i\in \{1, 2, \dots, n\}$, is
\[c(x_i(t)),\ i\in \{1, 2, \dots, n\}.\]
All firms have the same cost functions. It satisfies $c'(x_i(t))>0$. The instantaneous profit of Firm $i$, is
\[\pi_i(t)=x_i(t)p_i(x_1(t),x_2(t),\dots,x_n(t))-c(x_i(t)),\ i\in \{1, 2, \dots, n\}.\]
The moving of the number of firms is governed by
\begin{equation}
\frac{dn(t)}{dt}=s\left[\sum_{i=1}^nx_i(t)p_i(x_1(t),x_2(t),\dots,x_n(t))-\sum_{i=1}^nc(x_i(t))\right],\ s>0.\label{m1}
\end{equation}
The number of firms increases or decreases proportionally to the total profit of the firms.

The problem of Firm $i$ is
\[\max_{x_i(t)} \int_{0}^{\infty}e^{-\rho t}[x_i(t)p_i(x_1(t),x_2(t),\dots, x_n(t))-c(x_i(t))]dt,\]
subject to (\ref{m1}). $\rho>0$ is the discount rate.

The present value Hamiltonian function of Firm $i,\ i\in \{1, 2, \dots, n\}$, is
\begin{align*}
\mathcal{H}_i(t)=&e^{-\rho t}\left\{x_i(t)p(x_1(t),x_2(t),\dots, x_n(t))-c(x_i(t))\right.\\
&\left.+\lambda_i(t)s\left[\sum_{j=1}^nx_j(t)p_j(x_1(t),x_2(t),\dots,x_n(t))-\sum_{j=1}^nc(x_j(t))\right]\right\}.
\end{align*}
The current value Hamiltonian function of Firm $i,\ i\in \{1, 2, \dots, n\}$, is
\begin{align*}
\hat{\mathcal{H}}_i(t)=&e^{\rho t}\mathcal{H}_1(t)=x_i(t)p_i(x_1(t),x_2(t),\dots, x_n(t))-c(x_i(t))\\
&+\lambda_i(t)s\left[\sum_{j=1}^nx_j(t)p_j(x_1(t),x_2(t),\dots,x_n(t))-\sum_{j=1}^nc(x_j(t))\right].
\end{align*}
Let
\[\mu_i(t)=e^{-\rho t}\lambda_i(t),\ i\in \{1, 2, \dots, n\}.\]
$\mu_i(t)$ is the costate variable. 

Assume that the outputs of all firms are equal. The free entry condition is
\[p_i(x,x,\dots,x)x-c(x)-f=0.\]
From this
\begin{equation*}
\frac{dn}{dx}=-\frac{p_i(x,x,\dots,x)+\frac{\partial p_i}{\partial x_i(t)}x+(n-1)\frac{\partial p_i}{\partial x_j(t)}x-c'(x)}{\frac{\partial p_i}{\partial x_j(t)}x^2}.
\end{equation*}

Suppose that a monopolistic firm produce $n$ substitutable goods. It determines the output of each good. By symmetry we assume that the outputs of all goods are equal. Let $x$ be the output of each good. Its profit is $np_i(x, x, \dots, x)x-nc(x).$ The condition for profit maximization at $t$ in the static equilibrium is
\[n\left[p_i(x,x,\dots,x)+\frac{\partial p_i}{\partial x_i(t)}x+(n-1)\frac{\partial p_i}{\partial x_j(t)}x-c'(x)\right]=0.\]
If 
\[p_i(x,x,\dots,x)+\frac{\partial p_i}{\partial x_i(t)}x+(n-1)\frac{\partial p_i}{\partial x_j(t)}x-c'(x)\geq 0,\]
the output of each firm in the steady states of open-loop and closed-loop solutions should be smaller than (or equal to) the output of each good by the above monopolist. Therefore, we assume 
\[p_i(x,x,\dots,x)+\frac{\partial p_i}{\partial x_i(t)}x+(n-1)\frac{\partial p_i}{\partial x_j(t)}x-c'(x)<0.\]
Then, 
\begin{equation}
\frac{dn}{dx}<0.\label{nn}
\end{equation}
This holds in all cases.

We can assume
\[p_i-c'(x_i(t))>0,\ i\in \{1, 2, \dots, n\}.\]
This means that the price of the good is larger than the marginal cost of the firms. 

Consider a case such that each firm determines its output given the prices of the goods of other firms. Then, the profit maximization condition for Firm  $i$ in the static oligopoly is
\begin{equation}
p_i+\frac{\partial p_i}{\partial x_i(t)}x_i(t)+\sum_{j\neq i}\frac{\partial p_i}{\partial x_j(t)}x_i(t)\frac{dx_j(t)}{dx_i(t)}-c'(x_i(t))=0.\label{cc1}
\end{equation}
From the condition that $p_j(x_1, x_2,\dots, x_n)$ is constant for each $j\neq i$, we have
\[\frac{\partial p_j}{\partial x_j(t)}\frac{dx_j(t)}{dx_i(t)}+\sum_{k\neq i,j}\frac{\partial p_j}{\partial x_k(t)}\frac{dx_k(t)}{dx_i(t)}+\frac{\partial p_j}{\partial x_i(t)}=0.\]
By symmetry $\frac{\partial p_j}{\partial x_k(t)}=\frac{\partial p_j}{\partial x_i(t)}$ and $\frac{dx_k(t)}{dx_i(t)}=\frac{dx_j(t)}{dx_i(t)}$. Then,
\[\frac{dx_j(t)}{dx_i(t)}=-\frac{\frac{\partial p_j}{\partial x_i(t)}}{\frac{\partial p_j}{\partial x_j(t)}+(n-2)\frac{\partial p_j}{\partial x_i(t)}}.\]
Again by symmetry $\frac{\partial p_j}{\partial x_j(t)}=\frac{\partial p_i}{\partial x_i(t)}$, $\frac{\partial p_j}{\partial x_i(t)}=\frac{\partial p_i}{\partial x_j(t)}$ at the equilibrium. Thus, (\ref{cc1}) is rewritten as 
\begin{align*}
&p_i+\frac{\left(\frac{\partial p_i}{\partial x_i(t)}\right)^2+(n-2)\frac{\partial p_i}{\partial x_i(t)}\frac{\partial p_i}{\partial x_j(t)}-(n-1)\left(\frac{\partial p_i}{\partial x_j(t)}\right)^2}{\frac{\partial p_i}{\partial x_i(t)}+(n-2)\frac{\partial p_i}{\partial x_j(t)}}x_i(t)-c'(x_i(t))\\
=&p_i+\frac{\left(\frac{\partial p_i}{\partial x_i(t)}-\frac{\partial p_i}{\partial x_j(t)}\right)\left[\frac{\partial p_i}{\partial x_i(t)}+(n-1)\frac{\partial p_i}{\partial x_j(t)}\right]}{\frac{\partial p_i}{\partial x_i(t)}+(n-2)\frac{\partial p_i}{\partial x_j(t)}}x_i(t)-c'(x_i(t))=0.\notag
\end{align*}
Since $\frac{\partial p_i}{\partial x_i(t)}<0$, $\frac{\partial p_i}{\partial x_j(t)}<0$ and $\left|\frac{\partial p_i}{\partial x_i(t)}\right|>\left|\frac{\partial p_i}{\partial x_j(t)}\right|$, we have
\begin{equation}
\frac{\partial p_i}{\partial x_i(t)}-\frac{\partial p_i}{\partial x_j(t)}<0,\ \frac{\frac{\partial p_i}{\partial x_i(t)}+(n-1)\frac{\partial p_i}{\partial x_j(t)}}{\frac{\partial p_i}{\partial x_i(t)}+(n-2)\frac{\partial p_i}{\partial x_j(t)}}>1.\label{cc21}
\end{equation}
If 
\[p_i+\frac{\left(\frac{\partial p_i}{\partial x_i(t)}-\frac{\partial p_i}{\partial x_j(t)}\right)\left[\frac{\partial p_i}{\partial x_i(t)}+(n-1)\frac{\partial p_i}{\partial x_j(t)}\right]}{\frac{\partial p_i}{\partial x_i(t)}+(n-2)\frac{\partial p_i}{\partial x_j(t)}}x_i(t)-c'(x_i(t))\leq 0,\]
at the steady state of open-loop and closed-loop solutions, the output of each firm is larger than (or equal to) that under the above Bertrand type behaviors of firms. Thus, we assume
\[p_i+\frac{\left(\frac{\partial p_i}{\partial x_i(t)}-\frac{\partial p_i}{\partial x_j(t)}\right)\left[\frac{\partial p_i}{\partial x_i(t)}+(n-1)\frac{\partial p_i}{\partial x_j(t)}\right]}{\frac{\partial p_i}{\partial x_i(t)}+(n-2)\frac{\partial p_i}{\partial x_j(t)}}x_i(t)-c'(x_i(t))>0,\]
at the steady states of open-loop and closed-loop solutions of dynamic oligopoly. From (\ref{cc21}) we can assume\begin{equation}
p_i+\left(\frac{\partial p_i}{\partial x_i(t)}-\frac{\partial p_i}{\partial x_j(t)}\right)x_i(t)-c'(x_i(t))>0.\label{cc3}
\end{equation}

\section{The open-loop solution}

\subsection{General analysis}

We seek to the general open-loop solution. The first order condition for Firm $i$ is
\begin{align}
\frac{\partial \hat{\mathcal{H}}_i(t)}{\partial x_i(t)}=&p_i+\frac{\partial p_i}{\partial x_i(t)}x_i(t)-c'(x_i(t))\label{e6}\\
&+\lambda_{i}(t)s\left[p_i+\frac{\partial p_i}{\partial x_i(t)}x_i(t)-c'(x_i(t))+\sum_{j\neq i}\frac{\partial p_j}{\partial x_i(t)}x_j(t)\right]=0.\notag
\end{align}
The second order condition is
\begin{align}
\frac{\partial^2 \hat{\mathcal{H}}_i(t)}{\partial x_i(t)^2}=&2\frac{\partial p_i}{\partial x_i(t)}+\frac{\partial^2 p_i}{\partial x_i(t)^2}-c''(x_i(t))\label{s1}\\
&+\lambda_{i}(t)s\left[2\frac{\partial p_i}{\partial x_i(t)}+\frac{\partial^2 p_i}{\partial x_i(t)^2}x_i(t)-c''(x_i(t))+\sum_{j\neq i}\frac{\partial^2 p_j}{\partial x_i(t)^2}x_j(t)\right]<0.\notag
\end{align}
The adjoint condition is
\begin{align}
-&\frac{\partial \hat{\mathcal{H}}_i(t)}{\partial n(t)}=-\frac{\partial p_i}{\partial x_j(t)}x_i(t)x_j(t)\label{ad1}\\
&-\lambda_i(t)s\left[\frac{\partial p_j}{\partial x_k(t)}x_k(t)\sum_{l=1}^nx_l(t)+x_j(t)p_j-c(x_j(t)))\right]\notag\\
=&\frac{\partial \lambda_i(t)}{\partial t}-\rho\lambda_i(t),\ j\neq i,\ k\neq j.\notag
\end{align}
At the steady state we have $x_i(t)p\left(\sum_{j=1}^nx_j(t)\right)-c(x_i(t))=0$ and $\frac{\partial \lambda_i(t)}{\partial t}=0$ for all $i\in \{1,2, \dots, n\}$. By symmetry, all $x_i(t)$'s and all $\lambda_i(t)$'s are equal, and
\[\frac{\partial p_i}{\partial x_j(t)}=\frac{\partial p_j}{\partial x_i(t)}=\frac{\partial p_j}{\partial x_k(t)},\ j\neq i,\ k\neq j,\]
\[\frac{\partial^2 p_j}{\partial x_i(t)^2}=\frac{\partial^2 p_i}{\partial x_j(t)^2},\ j\neq i.\]
Denote the steady state values of $x_i(t)$, $\lambda_i(t)$ and $n(t)$ by $x^*$, $\lambda^*$ and $n^*$. (\ref{e6}), (\ref{s1}) and (\ref{ad1}) are reduced to
\begin{equation}
p_i+\frac{\partial p_i}{\partial x_i(t)}x^*-c'(x^*)+\lambda^*s\left(p_i+\frac{\partial p_i}{\partial x_i(t)}+(n^*-1)\frac{\partial p_i}{\partial x_j(t)}x^*-c'(x^*)\right)=0,\label{fo1}
\end{equation}
\[2\frac{\partial p_i}{\partial x_i(t)}+\frac{\partial^2 p_i}{\partial x_i(t)^2}x^*-c''(x^*)+\lambda^*s\left(2\frac{\partial p_i}{\partial x_i(t)}+\frac{\partial^2 p_i}{\partial x_i(t)^2}+(n^*-1)\frac{\partial^2 p_i}{\partial x_j(t)^2}x^*-c''(x^*)\right)<0,\]
and
\begin{align*}
-\frac{\partial p_i}{\partial x_j(t)}(x^*)^2-\lambda^*n^*s\frac{\partial p_i}{\partial x_j(t)}(x^*)^2=-\rho\lambda^*.
\end{align*}
Therefore,
\begin{align}
\lambda^*s=\frac{s\frac{\partial p_i}{\partial x_j(t)}(x^*)^2}{\rho-n^*s\frac{\partial p_i}{\partial x_j(t)}(x^*)^2}<0.\label{ad2}
\end{align}
From (\ref{fo1}) and (\ref{ad2})
\begin{align}
&p_i+\frac{\partial p_i}{\partial x_i(t)}x^*-c'(x^*) \label{fo12}\\
&+\frac{s\frac{\partial p_i}{\partial x_j(t)}(x^*)^2}{\rho-n^*s\frac{\partial p_i}{\partial x_j(t)}(x^*)^2}\left(p_i+\frac{\partial p_i}{\partial x_i(t)}x^*+(n^*-1)\frac{\partial p_i}{\partial x_j(t)}x^*-c'(x^*)\right)=0.\notag
\end{align}
Let $\tilde{x}$ and $\tilde{n}$ be the equilibrium output of each firm and the number of firms in the static monopolistic competition. Then,
\[p_i+\frac{\partial p_i}{\partial x_i(t)}\tilde{x}-c'(\tilde{x})=0.\]
Suppose that $x=\tilde{x}$ for each firm and $n=\tilde{n}$. The left-hand side of (\ref{fo12}) is
\[(n-1)\frac{s\frac{\partial p_i}{\partial x_j(t)}x^2}{\rho-ns\frac{\partial p_i}{\partial x_j(t)}x^2}\frac{\partial p_i}{\partial x_j(t)}x.\]
This is positive. Thus, under the assumption that the second order condition is satisfied, the output of each firm in the open-loop solution is larger than that at the static equilibrium, that is, $x^*>\tilde{x}$. 

From (\ref{nn}) $n^*<\tilde{n}$. We obtain the following result.
\begin{proposition}
The number of firms at the steady state in the open-loop solution of monopolistic competition is smaller than that at the static equilibrium of monopolistic competition.
\end{proposition}
Note that when $\rho\rightarrow +\infty$ or $s\rightarrow 0$, the steady state of open-loop solution approaches to the static equilibrium..

\subsection{A linear example}

Suppose that the inverse demand function for Firm $i$ is
\begin{equation*}
p_i(t)=a-x_i(t)-b\sum_{j\neq i}^nx_j(t).
\end{equation*}
$a$ is a positive constant, and $0<b<1$. Also suppose that the cost function of Firm $i,\ i\in \{1, 2, \dots, n\}$, is
\begin{equation*}
c(x_i(t))=cx_i(t)+f,\ c>0.
\end{equation*}
$f>0$ is the fixed cost. The moving of the number of firms is governed by
\begin{equation*}
\frac{dn(t)}{dt}=s\left[\left(a-x_i(t)-b\sum_{j\neq i}^nx_j(t)\right)\sum_{i=1}^nx_i(t)-c\sum_{i=1}^nx_i(t)-n(t)f\right],\ s>0.
\end{equation*}
The current value Hamiltonian function is
\begin{align*}
\hat{\mathcal{H}}_i(t)=&x_i(t)\left(a-x_i(t)-b\sum_{j\neq i}^nx_j(t)\right)-cx_i(t)-f\\
&+\lambda_i(t)s\left[\left(a-x_j(t)-b\sum_{k\neq j}^nx_k(t)\right)\sum_{j=1}^nx_j(t)-c\sum_{j=1}^nx_j(t)-n(t)f\right].
\end{align*}
The first order and the second order conditions for Firm $i,\ i\in \{1, 2, \dots, n\}$, are
\begin{align*}
\frac{\partial \hat{\mathcal{H}}_i(t)}{\partial x_i(t)}=\left(1+\lambda_i(t)s\right)\left(a-2x_i(t)-b\sum_{i\neq i}^nx_j(t)-c\right)-\lambda_i(t)sb\sum_{j\neq i}x_j(t)=0,
\end{align*}
and
\[\frac{\partial^2 \hat{\mathcal{H}}_i(t)}{\partial x_i(t)^2}=-2\left(1+\lambda_i(t)s\right)<0.\]
The adjoint condition for Firm $i,\ i\in \{1, 2, \dots, n\}$, is
\begin{align*}
-\frac{\partial \hat{\mathcal{H}}_i(t)}{\partial n(t)}=&bx_i(t)x_j(t)-\lambda_i(t)s\left[\left(a-x_j(t)-b\sum_{k\neq j}x_k(t)\right)x_j(t).\right.\\
&\left.-cx_j(t)-f-bx_k(t)\sum_{l=1}^nx_l(t)\right]=\frac{\partial \lambda_i(t)}{\partial t}-\rho\lambda_i(t).
\end{align*}
At the steady state we have $(a-\sum_{j=1}^nx_j(t))x_i(t)-cx_i(t)-f=0$ and $\frac{\partial \lambda_i(t)}{\partial t}=0$ for all $i\in \{1,2, \dots, n\}$. By symmetry, all $x_i(t)$'s and all $\lambda_i(t)$'s are equal. Denote the steady state values of $x_i(t)$, $\lambda_i(t)$ and $n(t)$ by $x^*$, $\lambda^*$ and $n^*$. Then, the above adjoint condition is reduced to
\[b(x^*)^2+\lambda^*n^*sb(x^*)^2=-\rho\lambda^*.\]
From this
\[\lambda^*=-\frac{b(x^*)^2}{\rho+n^*sb(x^*)^2}<0,\]
or
\[\lambda^*s=-\frac{bs(x^*)^2}{\rho+n^*sb(x^*)^2}<0.\]
Since $\rho>0$, when $s\rightarrow 0$ we have $\lambda^*s\rightarrow 0$. Similarly, when $\rho\rightarrow +\infty$ we have $\lambda^*s\rightarrow 0$.

\begin{figure}[ptb]
\centering
\includegraphics[height=8cm]{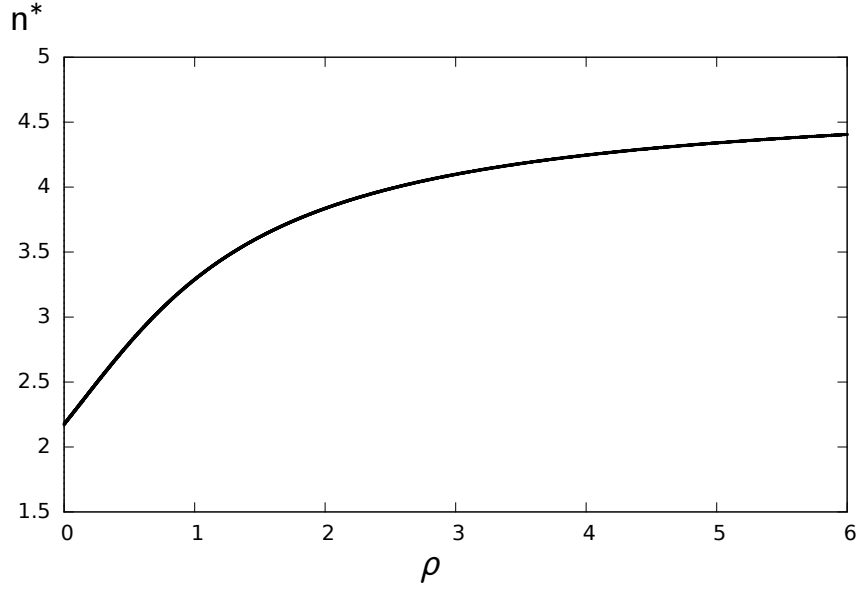}
	\vspace*{-.3cm}
	\caption{The numbers of firms in open-loop and $\rho$}
\end{figure}

At the steady state the first order condition is reduced to
\begin{align}
&\left(1+\lambda^*s\right)\left(a-2x^*-(n^*-1)bx^*-c\right)-\lambda^*s(n^*-1)bx^*\label{l1}\\
=&\left(a-2x^*-(n^*-1)bx^*-c\right)-\frac{b^2s(x^*)^2}{\rho+n^*bs(x^*)^2}\left(a-2^*x^*-2(n^*-1)bx^*-c\right)=0.\notag
\end{align}
On the other hand, the free entry condition at the steady state is
\begin{equation}
(a-x^*-(n^*-1)bx^*)x^*-c(x^*)-f=0.\label{l2}
\end{equation}
Solving (\ref{l1}) and (\ref{l2}) we get the steady state values of $x^*$ and $n^*$. We give graphical representations in Figure 1 assuming $a=11,\ f=4,\ c=1,\ s=\frac{1}{10},\ b=\frac{4}{5}$ and in Figure 2 assuming $a=11,\ f=4,\ c=1,\ \rho=\frac{1}{2},\ b=\frac{4}{5}$.

\begin{figure}[ptb]
\centering
\includegraphics[height=8cm]{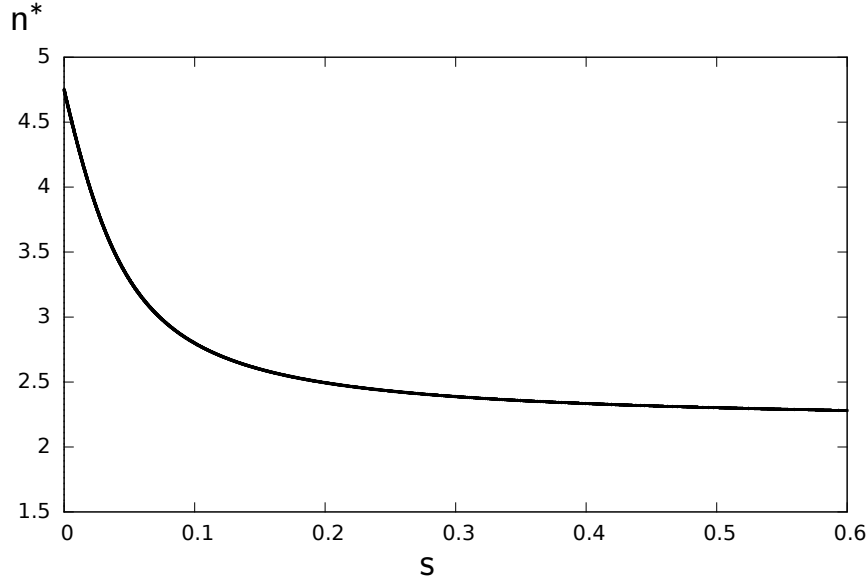}
	\vspace*{-.3cm}
	\caption{The numbers of firms in open-loop and $s$}
\end{figure}

When $s\rightarrow 0$ or $\rho\rightarrow +\infty$, (\ref{l1}) is further reduced to
\[a-2x^*-(n^*-1)bx^*-c=0.\]
This is equivalent to the static equilibrium condition.

\section{The memoryless closed-loop solution: A general analysis}

We seek to a memoryless closed-loop solution. The first order condition and the second order condition are the same as those in the open-loop solution as follows.
\begin{align}
\frac{\partial \hat{\mathcal{H}}_i(t)}{\partial x_i(t)}=&p_i+\frac{\partial p_i}{\partial x_i(t)}x_i(t)-c'(x_i(t))\label{e7}\\
&+\lambda_{i}(t)s\left[p_i+\frac{\partial p_i}{\partial x_i(t)}x_i(t)-c'(x_i(t))+\sum_{j\neq i}\frac{\partial p_j}{\partial x_i(t)}x_j(t)\right]=0,\notag
\end{align}
and
\begin{align*}
\frac{\partial^2 \hat{\mathcal{H}}_i(t)}{\partial x_i(t)^2}=&2\frac{\partial p_i}{\partial x_i(t)}+\frac{\partial^2 p_i}{\partial x_i(t)^2}-c''(x_i(t))\\
&+\lambda_{i}(t)s\left[2\frac{\partial p_i}{\partial x_i(t)}+\frac{\partial^2 p_i}{\partial x_i(t)^2}x_i(t)-c''(x_i(t))+\sum_{j\neq i}\frac{\partial^2 p_j}{\partial x_i(t)^2}x_j(t)\right]<0.\notag
\end{align*}
The adjoint condition is different from that in the open-loop solution. It is written as
\begin{align}
-\frac{\partial \hat{\mathcal{H}}_i(t)}{\partial n(t)}-\sum_{j\neq i}\frac{\partial \hat{\mathcal{H}}_i(t)}{\partial x_j(t)}\frac{\partial x_i(t)}{\partial n(t)}=\frac{\partial \lambda_i(t)}{\partial t}-\rho\lambda_i(t).\label{ad3}
\end{align}
The term in (\ref{ad3})
\[-\sum_{j\neq i}\frac{\partial \hat{\mathcal{H}}_i(t)}{\partial x_j(t)}\frac{\partial x_i(t)}{\partial n(t)}\]
takes into account the interaction between the control variable of the firms other than Firm $i$ and the current level of the state variable. We have
\begin{align*}
\sum_{j\neq i}\frac{\partial \hat{\mathcal{H}}_i(t)}{\partial x_j(t)}=&(n-1)x_i(t)\frac{\partial p_i}{\partial x_j(t)}\\
&+\lambda_i(t)s\sum_{j\neq i}\left[p_j+\frac{\partial p_j}{\partial x_j(t)}x_j(t)+\sum_{k\neq j}^n\frac{\partial p_k}{\partial x_j(t)}x_k(t)-c'(x_j(t))\right], j\neq i,\ k\neq j.
\end{align*}
From (\ref{e7})
\begin{align}
\frac{\partial x_i(t)}{\partial n(t)}=&-\frac{1}{\Delta}\left\{(1+\lambda_i(t)s)\left(\frac{\partial p_i}{\partial x_j(t)}+\frac{\partial^2 p_i}{\partial x_i(t)\partial x_j(t)}x_i(t)\right)x_j(t)\right.\label{xi}\\
&\left.+\lambda_i(t)s\left(\frac{\partial p_i}{\partial x_j(t)}x_j(t)+\sum_{j\neq i}\frac{\partial p_j}{\partial x_i(t)\partial x_k(t)}x_j(t)x_k(t)\right)\right\},\ j\neq i,\ k\neq i,j,\notag
\end{align}
where
\begin{align}
\Delta=&2\frac{\partial p_i}{\partial x_i(t)}+\frac{\partial^2 p_i}{\partial x_i(t)^2}-c''(x_i(t))\label{del}\\
&+\lambda_{i}(t)s\left[2\frac{\partial p_i}{\partial x_i(t)}+\frac{\partial^2 p_i}{\partial x_i(t)^2}x_i(t)-c''(x_i(t))+\sum_{j\neq i}\frac{\partial^2 p_j}{\partial x_i(t)^2}x_j(t)\right]<0.\notag
\end{align}
At the steady state we have $p\left(\sum_{k=1}^nx_k(t)\right)x_i(t)-c(x_i(t))=0$ and $\frac{\partial \lambda_i(t)}{\partial t}=0$ for all $i\in \{1,2, \dots, n\}$. By symmetry, all $x_i(t)$'s and all $\lambda_i(t)$'s are equal. Denote the steady state values of $x_i(t)$, $\lambda_i(t)$ and $n(t)$ by $x^{**}$, $\lambda^{**}$ and $n^{**}$. Then, using $\frac{\partial p_j}{\partial x_i(t)}=\frac{\partial p_i}{\partial x_j(t)}$, (\ref{e7}) and (\ref{ad3}) are reduced to
\begin{align}
&p_i+\frac{\partial p_i}{\partial x_i(t)}x^{**}-c'(x^{**})\label{cl31}\\
&+\lambda^{**}s\left[p_i+\frac{\partial p_i}{\partial x_i(t)}x^{**}+(n^{**}-1)\frac{\partial p_i}{\partial x_j(t)}x^{**}-c'(x^{**})\right]=0,\notag
\end{align}
and
\begin{align}
&-\frac{\partial p_i}{\partial x_j(t)}(x^{**})^2-\lambda^{**}n^{**}s\frac{\partial p_i}{\partial x_j(t)}(x^{**})^2-\left[(n^{**}-1)\frac{\partial p_i}{\partial x_j(t)}x^{**}\right.\label{cl4}\\
&\left.+\lambda^{**}s(n^{**}-1)\left(p_i+\frac{\partial p_i}{\partial x_i(t)}x^{**}+(n^{**}-1)\frac{\partial p_i}{\partial x_j(t)}x^{**}-c'(x^{**})\right)\right]\frac{\partial x_i(t)}{\partial n(t)}\notag\\
=&-\rho\lambda^{**}.\notag
\end{align}
From (\ref{cl31})
\begin{equation}
\lambda^{**}s=-\frac{p_i+\frac{\partial p_i}{\partial x_i(t)}x^{**}-c'(x^{**})}{p_i+\frac{\partial p_i}{\partial x_i(t)}x^{**}+(n^{**}-1)\frac{\partial p_i}{\partial x_j(t)}x^{**}-c'(x^{**})},\label{e71}
\end{equation}
and
\begin{equation}
1+\lambda^{**}s=\frac{(n^{**}-1)\frac{\partial p_i}{\partial x_j(t)}x^{**}}{p_i+\frac{\partial p_i}{\partial x_i(t)}x^{**}+(n^{**}-1)\frac{\partial p_i}{\partial x_j(t)}x^{**}-c'(x^{**})}.\label{e72}
\end{equation}
Then, (\ref{cl4}) is rewritten as
\begin{align*}
&-\frac{\partial p_i}{\partial x_j(t)}(x^{**})^2-\lambda^{**}n^{**}s\frac{\partial p_i}{\partial x_j(t)}(x^{**})^2\\
&+(n^{**}-1)\left[p_i+\left(\frac{\partial p_i}{\partial x_i(t)}-\frac{\partial p_i}{\partial x_j(t)}\right)x^{**}-c'(x^{**})\right]\frac{\partial x_i(t)}{\partial n(t)}=-\rho\lambda^{**}.
\end{align*}
This means
\begin{equation}
\lambda^{**}s=\frac{\frac{\partial p_i}{\partial x_j(t)}s(x^{**})^2-(n^{**}-1)s\left[p_i+\left(\frac{\partial p_i}{\partial x_i(t)}-\frac{\partial p_i}{\partial x_j(t)}\right)x^{**}-c'(x^{**})\right]\frac{\partial x_i(t)}{\partial n(t)}}{\rho-n^{**}s\frac{\partial p_i}{\partial x_j(t)}(x^{**})^2}.\label{fo13}
\end{equation}
By (\ref{cc3}),
\[p_i+\left(\frac{\partial p_i}{\partial x_i(t)}-\frac{\partial p_i}{\partial x_j(t)}\right)x^{**}-c'(x^{**})>0.\]
From (\ref{fo13}) and (\ref{cl31}) we get
\begin{align}
&\frac{\partial \hat{\mathcal{H}}_i}{\partial x_i(t)}=p_i+\frac{\partial p_i}{\partial x_i(t)}x^{**}-c'(x^{**})\label{cl61}\\
&+\frac{\frac{\partial p_i}{\partial x_j(t)}s(x^{**})^2-(n^{**}-1)s\left[p_i+\left(\frac{\partial p_i}{\partial x_i(t)}-\frac{\partial p_i}{\partial x_j(t)}\right)x^{**}-c'(x^{**})\right]\frac{\partial x_i(t)}{\partial n(t)}}{\rho-n^{**}s\frac{\partial p_i}{\partial x_j(t)}(x^{**})^2}\left[p_i+\frac{\partial p_i}{\partial x_i(t)}x^{**}\right.\notag\\
&\left.+(n^{**}-1)\frac{\partial p_i}{\partial x_j(t)}x^{**}-c'(x^{**})\right]=0.\notag
\end{align}
Compare (\ref{cl61}) and (\ref{fo12}). Suppose that $x=x^{**}$, $n=n^{**}$ and (\ref{cl61}) is satisfied, the left-hand side of (\ref{fo12}) is
\begin{align}
&\frac{(n^{**}-1)s\left[p_i+\left(\frac{\partial p_i}{\partial x_i(t)}-\frac{\partial p_i}{\partial x_j(t)}\right)x^{**}-c'(x^{**})\right]\frac{\partial x_i(t)}{\partial n(t)}}{\rho-n^{**}s\frac{\partial p_i}{\partial x_j(t)}(x^{**})^2}\left[p_i+\frac{\partial p_i}{\partial x_i(t)}x^{**}\right.\label{e75}\\
&\left.+(n^{**}-1)\frac{\partial p_i}{\partial x_j(t)}x^{**}-c'(x^{**})\right].\notag
\end{align}
At the steady state from (\ref{xi})
\begin{align}
&\frac{\partial x_i(t)}{\partial n(t)}=\frac{1}{\Gamma}\left\{-(n^{**}-1)\frac{\partial p_i}{\partial x_j(t)}x^{**}\left(\frac{\partial p_i}{\partial x_j(t)}+\frac{\partial^2 p_i}{\partial x_i(t)\partial x_j(t)}x^{**}\right)x^{**}\right.\label{e74}\\
&\left.+\left(p_i+\frac{\partial p_i}{\partial x_i(t)}x^{**}-c'(x^{**})\right)\left[\frac{\partial p_i}{\partial x_j(t)}+(n^{**}-1)\frac{\partial^2 p_i}{\partial x_j(t)\partial x_k(t)}x^{**}\right]x^{**}\right\},\notag\\
=&\frac{1}{\Gamma}\left\{(n^{**}-1)(p_i-c'(x^{**}))\left(\frac{\partial p_i}{\partial x_j(t)}+\frac{\partial^2 p_i}{\partial x_j(t)\partial x_k(t)}x^{**}\right)x^{**}\right.\notag\\
&\left.-(n^{**}-1)\frac{\partial p_i}{\partial x_j(t)}\left[\left(\frac{\partial p_i}{\partial x_j(t)}+\frac{\partial^2 p_i}{\partial x_i(t)\partial x_j(t)}x^{**}\right)-\left(\frac{\partial p_i}{\partial x_j(t)}+\frac{\partial^2 p_i}{\partial x_i(t)\partial x_j(t)}x^{**}\right)\right]x^{**}\right.\notag\\
&\left.-(n^{**}-2)\frac{\partial p_i}{\partial x_j(t)}\left(p_i+\frac{\partial p_i}{\partial x_i(t)}x^{**}-c'(x^{**})\right)x^{**}\right\},\notag
\end{align}
where
\[\Gamma=\Delta\left(p_i+\frac{\partial p_i}{\partial x_i(t)}x^{**}+(n^{**}-1)\frac{\partial p_i}{\partial x_j(t)}x^{**}-c'(x^{**})\right).\]
We have $\Delta<0$, $p_i+\frac{\partial p_i}{\partial x_i(t)}x^{**}+(n^{**}-1)\frac{\partial p_i}{\partial x_j(t)}x^{**}-c'(x^{**})<0$, $p_i-c'(x^{**})>0$, and from (\ref{sub1}), (\ref{sub2}) and (\ref{sub3}), \[\frac{\partial p_i}{\partial x_j(t)}+\frac{\partial^2 p_i}{\partial x_i(t)\partial x_j(t)}x_i(t)<0,\]
\[\frac{\partial p_i}{\partial x_j(t)}+\frac{\partial^2 p_i}{\partial x_j(t)\partial x_k(t)}x_i(t)<0,\]
and
\[\left|\frac{\partial p_i}{\partial x_j(t)}+\frac{\partial^2 p_i}{\partial x_i(t)\partial x_j(t)}x_i(t)\right|\geq \left|\frac{\partial p_i}{\partial x_j(t)}+\frac{\partial^2 p_i}{\partial x_j(t)\partial x_k(t)}x_i(t)\right|.\]
Suppose $\frac{\partial x_i(t)}{\partial n(t)}\geq 0$. From (\ref{fo13}), $\lambda^*s\leq 0$. By (\ref{cl31}), we have $p_i+\frac{\partial p_i}{\partial x_i(t)}x^{**}-c'(x^{**})\leq 0$.  From (\ref{e74}) this means $\frac{\partial x_i(t)}{\partial n(t)}<0$. It is a contradiction. Thus, we have $\frac{\partial x_i(t)}{\partial n(t)}<0$, and then (\ref{e75}) is positive (because $p_i+\frac{\partial p_i}{\partial x_i(t)}x^{**}+(n^{**}-1)\frac{\partial p_i}{\partial x_j(t)}x^{**}-c'(x^{**})<0$). This means $x^{**}<x^{*}$ and $n^{**}>n^*$. We have shown the following result.
\begin{proposition}
The number of firms at the steady state  in the memoryless closed-loop solution of monopolistic competition is larger than that in the open-loop solution of monopolistic competition.
\end{proposition}

If $\frac{\partial p_i}{\partial x_j(t)}s(x^{**})^2-(n^{**}-1)s\left[p_i+\left(\frac{\partial p_i}{\partial x_i(t)}-\frac{\partial p_i}{\partial x_j(t)}\right)x^{**}-c'(x^{**})\right]\frac{\partial x_i(t)}{\partial n(t)}>0$ in (\ref{cl61}), $x^{**}<\tilde{x}$ and the number of firms at the steady state in the closed-loop solution is larger than that at the static equilibrium of free entry oligopoly.

Also note that from (\ref{cl61}) we find that when $s\rightarrow 0$ or $\rho\rightarrow +\infty$, the steady state of the closed-loop solution approaches to the static equilibrium.

\section{The memoryless closed-loop solution: A linear example}

\begin{figure}[ptb]
\centering
\includegraphics[height=8cm]{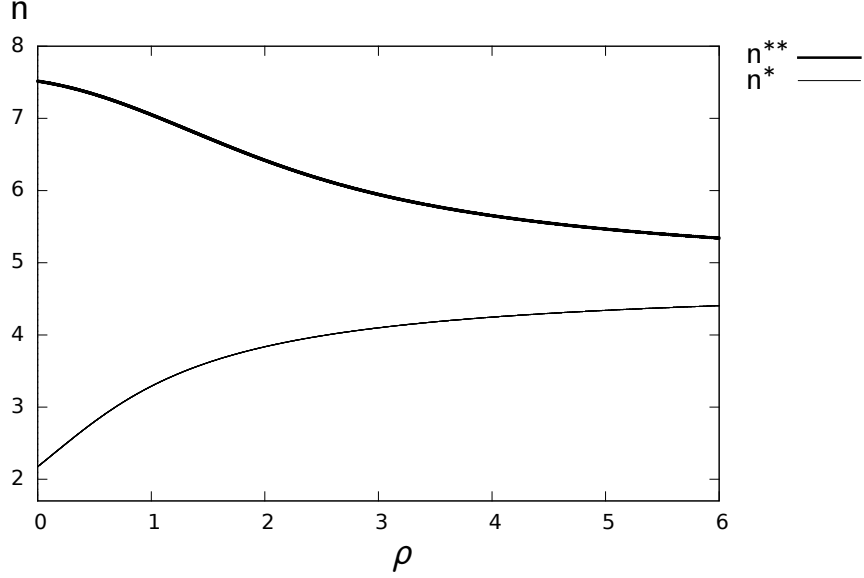}
	\vspace*{-.3cm}
	\caption{The numbers of firms in open-loop, closed-loop and $\rho$}
\end{figure}

Similarly to the example in the open-loop case, we assume that the inverse demand function is
\begin{equation*}
p_i(t)=a-x_i(t)-b\sum_{j\neq i}^nx_j(t),\ a>0,
\end{equation*}
and the cost function of Firm $i,\ i\in \{1, 2, \dots, n\}$, is
\begin{equation*}
c(x_i(t))=cx_i(t)+f,\ c>0,\ f>0.
\end{equation*}
The moving of the number of firms is governed by
\begin{equation*}
\frac{dn(t)}{dt}=s\left[\left(a-x_i(t)-b\sum_{j\neq i}^nx_j(t)\right)\sum_{i=1}^nx_i(t)-c\sum_{i=1}^nx_i(t)-n(t)f\right],\ s>0.
\end{equation*}
From (\ref{e71}), (\ref{e72}) and (\ref{del}), at the steady state we have
\[\lambda^{**}s=-\frac{a-2x^{**}-(n^{**}-1)bx^{**}-c}{a-2x^{**}-2(n^{**}-1)bx^{**}-c},\]
\[1+\lambda^{**}s=-\frac{(n^{**}-1)bx^{**}}{a-2x^{**}-2(n^{**}-1)bx^{**}-c},\]
\[\Delta=\frac{2(n^{**}-1)x^{**}}{a-2x^{**}-2(n^{**}-1)bx^{**}-c},\]
and
\[\Delta[a-2x^{**}-2(n^{**}-1)bx^{**}-c)]=2(n^{**}-1)x^{**}.\]
Therefore,
\begin{align*}
\frac{\partial x_i(t)}{\partial n(t)}=&-\frac{(n^{**}-1)(a-x^{**}-(n^{**}-1)bx^{**}-c)b-(n^{**}-2)(a-2x^{**}-(n^{**}-1)bx^{**}-c)b}{2(n^{**}-1)}\\
=&-\frac{[a+(n^{**}-3)x^{**}-(n^{**}-1)bx^{**}-c]b}{2(n^{**}-1)},
\end{align*}
(\ref{cl61}) is reduced to
\begin{align}
&\left(a-2x^{**}-(n^{**}-1)bx^{**}-c\right)\label{cl11}\\
&-\frac{sb(x^{**})^2-s(a-n^{**}bx^{**}-c)\frac{[a+(n^{**}-3)x^{**}-(n^{**}-1)bx^{**}-c]b}{2}}{\rho+n^{**}bs(x^{**})^2}\left(a-2n^{**}x^{**}-c\right)=0.\notag
\end{align}
On the other hand, the free entry condition at the steady state is the same as that in the open-loop case as follows,
\begin{equation}
(a-n^{**}-x^{**})x^{**}-c(x^{**})-f=0.\label{cl21}
\end{equation}
Solving (\ref{cl11}) and (\ref{cl21}) we get the steady state values of $x^{**}$ and $n^{**}$. We give graphical representations in Figure 3 assuming $a=11,\ f=4,\ c=1,\ s=\frac{1}{10},\ b=\frac{4}{5}$ and in Figure 4 assuming $a=11,\ f=4,\ c=1,\ \rho=\frac{1}{2},\ b=\frac{4}{5}$. In these figures we depict the relations between the number of firms at the steady states of open-loop and closed-loop solutions and the value of $s$ or $\rho$.

\begin{figure}[ptb]
\centering
\includegraphics[height=8cm]{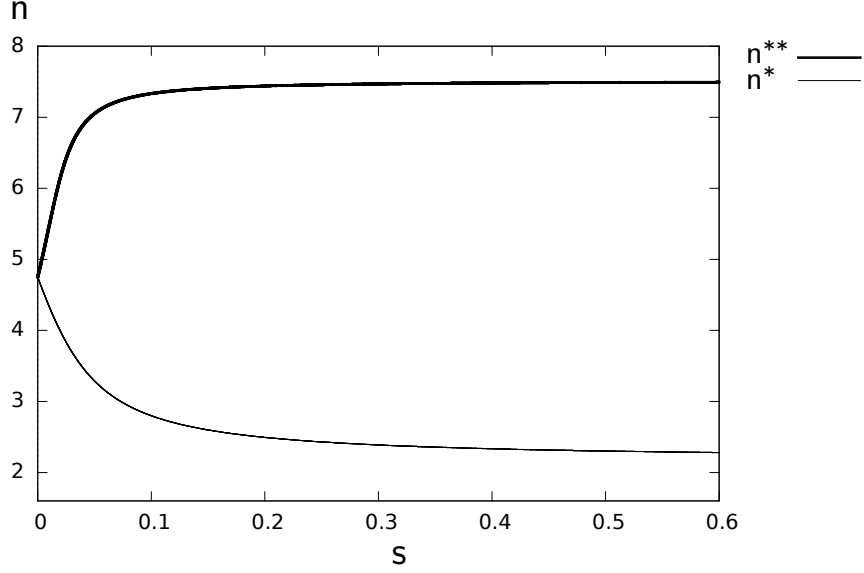}
	\vspace*{-.3cm}
	\caption{The numbers of firms in open-loop, closed-loop and $s$}
\end{figure}

\section{Concluding Remark}

In this paper we analyze a dynamic free entry oligopoly with differentiated goods, that is, a monopolistic competition by differential game approach.

\section*{Acknowledgment}

This work was supported by Japan Society for the Promotion of Science KAKENHI Grant Number 18K01594.

\bibliographystyle{plainnat} 
\bibliography{yatanaka}

\newpage

\end{document}